\documentclass[12pt]{amsart}

\usepackage{amsmath}
\usepackage{graphicx}
\usepackage{latexsym}
\usepackage[dvipsnames]{xcolor}
\usepackage{amscd}
\usepackage[all]{xy}
\usepackage{enumerate}
\usepackage{hyperref}
\usepackage{subfigure}
\usepackage{soul}
\usepackage{comment}
\usepackage{epsfig}
\usepackage{epstopdf}
\usepackage{tikz}
\usepackage{marginnote}
\usepackage{lipsum}
\usepackage{mathtools}
\usepackage{multirow}

\newtheorem{thm}{Theorem}

\newtheorem{theorem}[thm]{Theorem}

\newtheorem{proposition}[thm]{Proposition}

\theoremstyle{definition}

\newtheorem*{definition*}{Definition}

\newtheorem{definition}[thm]{Definition}

\newtheorem{rk}[thm]{Remark}

\makeatletter
\newsavebox\myboxA
\newsavebox\myboxB
\newlength\mylenA

\newcommand*\xoverline[2][0.75]{%
    \sbox{\myboxA}{$\m@th#2$}%
    \setbox\myboxB\null
    \ht\myboxB=\ht\myboxA%
    \dp\myboxB=\dp\myboxA%
    \wd\myboxB=#1\wd\myboxA
    \sbox\myboxB{$\m@th\overline{\copy\myboxB}$}
    \setlength\mylenA{\the\wd\myboxA}
    \addtolength\mylenA{-\the\wd\myboxB}%
    \ifdim\wd\myboxB<\wd\myboxA%
       \rlap{\hskip 0.5\mylenA\usebox\myboxB}{\usebox\myboxA}%
    \else
        \hskip -0.5\mylenA\rlap{\usebox\myboxA}{\hskip 0.5\mylenA\usebox\myboxB}%
    \fi}
\makeatother

\newcommand{\CP}{{\mathbb{CP}}{}^2}
\newcommand{\CPb}{{\xoverline[0.75]{\mathbb{CP}}}{}^2}
\newcommand{\CPo}{{\mathbb{CP}}^{1}}
\newcommand{\RP}{{\mathbb{RP}}^{2}}
\newcommand{\R}{\mathbb{R}}

\newcommand{\N}{\mathbb{N}}
\newcommand{\Z}{\mathbb{Z}}

\newcommand{\K}{{\mathrm{K}3}}
\newcommand{\Kb}{\xoverline[0.75]{\mkern1mu{\mathrm{K}3}}\mkern-1mu}
\newcommand{\Xb}{\xoverline[0.7]{\mkern-1.1mu{X}\mkern1.1mu}}

\newcommand{\id}{\mathop{\mathrm{id}}\nolimits}
 
\newcommand{\Int}{\mathop{\mathrm{Int}}\nolimits} 
\newcommand{\Bd}{\partial}

\newcommand{\cs}{\mathbin{\#}}
\newcommand{\csb}{\mathbin{\sharp}}

\renewcommand{\:}{\,{:}\;}

\def \x {\times}

\def\R{{\mathbb{R}}}
\def\Z{{\mathbb{Z}}}

\def\N{{\mathbb{N}}}

\def\C{{\mathbb{C}}}

\begin{document}

\title[Branched covering $4$-manifolds]{Branched covering simply-connected $4$-manifolds}

\author[D. Auckly]{David Auckly} 
\address{Mathematics Department, Kansas State University, Manhattan, KS 66506-2602, USA}
\email{dav@ksu.edu}

\author[R. \.{I}. Baykur]{R. \.{I}nan\c{c} Baykur}
\address{Department of Mathematics and Statistics, University of Massachusetts, Amherst, MA 01003-9305, USA}
\email{baykur@math.umass.edu}

\author[R. Casals]{Roger Casals}
\address{University of California Davis, Dept. of Mathematics, Shields Avenue, Davis, CA 95616, USA}
\email{casals@math.ucdavis.edu}

\author[S. Kolay]{\\ Sudipta Kolay}
\address{School of Mathematics, Georgia Institute of Technology, Atlanta, GA 30332, USA}
\email{skolay3@math.gatech.edu}

\author[T. Lidman]{Tye Lidman}
\address{Department of Mathematics, North Carolina State University, Raleigh, NC, 27607, USA}
\email{tlid@math.ncsu.edu}

\author[D. Zuddas]{Daniele Zuddas}
\address{Dipartimento di Matematica e Geoscienze, Universit\`a di Trieste, Via Valerio 12/1, 34127 Trieste, Italy}
\email{dzuddas@units.it}

\begin{abstract} 
We prove that any closed simply-connected smooth $4$-manifold is $16$-fold branched covered by a product of an orientable surface with the $2$-torus,  where the construction is natural with respect to spin structures.  In particular this solves Problem 4.113(C) in Kirby's list. We also discuss analogous results for other families of $4$-manifolds with infinite fundamental groups. 
\end{abstract}


\maketitle

\setcounter{secnumdepth}{2}
\setcounter{section}{0}
\vspace{0.5in}

\section{Introduction} 

The work in this note was prompted by the following natural question:

\noindent \textit{``Does every closed $4$--manifold admit a branched covering by a symplectic $4$--manifold?'',}

\noindent which was studied by the authors at the 2018 American Institute of Mathematics Workshop on ``Symplectic four-manifolds through branched coverings'', motivated by a conjecture of Eliashberg \cite[Conjecture 6.2]{Eliashberg}. We provide a fairly strong answer to the above question in the case of simply connected $4$--manifolds,  which in particular solves Problem 4.113(C) in Kirby's list~\cite{KirbyList}:

\begin{theorem} \label{thm1}
Let $X$ be a closed oriented simply-connected smooth $4$--manifold. Then there exists $g\in\N$ and a degree $16$ branched covering $f\: X' \to X$ such that $X'$ is the smooth 4-manifold $T^2 \x \Sigma_g$.
In addition, if the 4--manifold $X$ is spin, the branched covering $f$ is natural with respect to a spin structure on $T^2 \x \Sigma_g$.
\end{theorem}

\noindent Note that the smooth 4-manifold $T^2 \x \Sigma_g$ admits a symplectic structure. It follows from Theorem \ref{thm1} that if instead $X$ is a closed (possibly non-orientable) connected smooth \linebreak $4$--manifold with finite $\pi_1(X)$, then there is a branched covering  $T^2 \x \Sigma_g \to X$ of degree $16 \, |\pi_1(X)|$, which factors through the universal covering $\widetilde{X} \to X$. 

The above do not generalize to $4$--manifolds with arbitrary fundamental groups; for instance, no $\Sigma_g$--bundle over $\Sigma_h$ with $g, h \geq 2$ and infinite monodromy group (e.g. any surface bundle with nonzero signature) can be dominated by a product $4$--manifold  by \mbox{\cite[Theorem~1.4]{KotschickLoh}.} Nonetheless, there are comparable results for many other $4$--manifolds with infinite fundamental groups. For example,  $X=\cs_g (S^1 \x S^3)$, with $\pi_1(X)\cong *\,\Z^{\,g}$, is degree $4$ branched covered by $X'=S^2 \x \Sigma_g$ by  \cite[Theorem 1.2]{PiergalliniZuddas}. In addition, the branched virtual fibering theorem of Sakuma \cite[Addendum~1]{Sakuma} implies the following:

\begin{proposition}  \label{prop1} 
Let $X=S^1 \x Y$ be a smooth \mbox{$4$--manifold} which is the product of $S^1$ and a closed connected oriented $3$--manifold $Y$. Then there exists $g\in\N$ and a double branched covering $X' \to X$, where $X'$ is a symplectic $4$--manifold which is a $\Sigma_g$--bundle over $T^2$.
\end{proposition}
\noindent Indeed, \cite{Sakuma} shows that any closed oriented \mbox{$3$--manifold} $Y$ is double branched covered by a surface bundle over a circle (see also \cite{Montesinos} for a different proof), from which Proposition \ref{prop1} is immediately deduced; this provides yet another class of $4$--manifolds with infinite fundamental group for which a (symplectic) branched cover can be readily described.  Here, we recall that the product of a fibered $3$--manifold and the circle is symplectic \cite{Thurston}.  

It is worth noting that with a little more information on the smooth topology of $X$, one can easily determine the topology of the branched coverings $X' \to X$ in Theorem~\ref{thm1} and Proposition~\ref{prop1}. For the former, one only needs to know the number of stabilizations by taking the connected sum with $S^2 \x S^2$ that are required before the simply-connected $4$--manifold $X$ completely decomposes into a connected sum of copies of $\CP$, $S^2 \x S^2$ and the $\K$ surface, taken with either orientation. This of course can always be achieved by a classical result of Wall \cite{Wall}, and for vast families of simply-connected $4$--manifolds, one stabilization is known to be enough \cite{BaykurSunukjian}. Similarly, for Proposition~\ref{prop1},  one just needs to know a Heegaard decomposition of the $3$--manifold factor $Y$ \cite{Sakuma}, or any open book on it \cite{Montesinos}. See Remark \ref{explicit} for some explicit examples. 

In all the results we have discussed above, the covering symplectic \mbox{$4$--manifold} $X'$  is not of general type, in contrast with the symplectic domination results of Fine--Panov \cite{FinePanov}; see Remark~\ref{domination} below. It would be  interesting to find more general families of non-symplectic \mbox{$4$--manifolds} branched covered (with universally fixed degree) by  specific families of symplectic $4$--manifolds like ours, say by $\Sigma_g$--bundles over $\Sigma_h$, for arbitrary $h$.

\vspace{0.2in}
\noindent \textit{Acknowledgments.} This project was started at the 2018 AIM workshop on \textit{``Symplectic four-manifolds through branched coverings'',} and was resumed following the 2020 BIRS Workshop on  \textit{``Interactions of gauge theory with contact and symplectic topology in dimensions 3 and 4."} The authors would like to thank the American Institute of Mathematics and the Banff International Research Station, and the other organizers of these workshops. D.\,A. was partially supported by the Simons Foundation grant 585139 and NSF grant DMS 1952755.  R.\,I.\,B. was  partially supported by the NSF grants DMS-200532 and DMS-1510395. R.\,C. is supported by the NSF grant DMS-1841913, the NSF CAREER grant DMS-1942363 and the Alfred P. Sloan Foundation. T.\,L. was partially supported by the NSF grant DMS-1709702 and a Sloan Fellowship. D.\,Z. was partially supported by the 2013 ERC Advanced Research Grant 340258 TADMICAMT; he is member of GNSAGA, Istituto Nazionale di Alta Matematica ``Francesco Severi'', Italy.
We would like to thank the anonymous referee for helpful comments.

\medskip
\section{Proof of Theorem \ref{thm1}}

Henceforth all the manifolds and maps we consider are assumed to be smooth. We denote by $\Xb$  the oriented $4$--manifold $X$ with the reversed orientation, and by $\cs_a X \cs_b Y$ the smooth connected sum of $a$ copies of $X$ and $b$ copies of $Y$. We denote by $\Sigma_g^b$ a closed connected oriented surface of genus $g$ with $b$ boundary components, and we drop $b$ from the notation when there is no boundary.

\subsection{Preliminaries}
Let us briefly recall the definition of a branched covering.
\begin{definition}
Let $X$ and $X'$ be compact connected smooth manifolds (possibly with boundary) of the same dimension, and let $f \: X' \to X$ be a smooth proper surjective map. We say that $f$ is a branched covering if it is finite-to-one and open, and moreover the (open) subset of $X'$ where $f$ is locally injective coincides with the subset of $X'$ where $f$ is a local diffeomorphism.\hfill$\Box$
\end{definition}
The subset $B'_f \subset X'$ where $f$ fails to be locally injective is called the \textit{branch set} of $f$, and its image $B_f = f(B'_f) \subset X$ is called the \textit{branch locus} of $f$. By a result of Church \cite[Corollary 2.3]{Church}, either $B'_f = \emptyset$ or $\dim B'_f = \dim B_f = \dim X -2$, and then the restriction of $f$ over the complement of $B_f$ is an ordinary connected covering space $X' \setminus f^{-1}(B_f) \to X \setminus B_f$.

Moroever, for every \textit{smooth} point of $B'_f$ at which $f|_{B'_f} \: B'_f \to X$ is a local smooth embedding, the map $f$ is \textit{topologically} locally equivalent to the map $p_d \: \C \times \R^{n-2} \to \C \times \R^{n-2}$ defined by $p_d (z, x) = (z^d, x)$, for some $d \geq 2$, where $n = \dim X' = \dim X$. However, the branched coverings $f_i$ that we consider below turn out to be \textit{smoothly} locally equivalent to $p_2$, while their composition, which will be indicated by $f$, has this property away from the singular points of $B'_f$. Notice that every finite composition of branched coverings is a branched covering, and the restriction to the boundary of a branched covering is a branched covering as well. Throughout, we assume that branched coverings between oriented manifolds are orientation-preserving.

\subsection{The argument} Let $X$ be a closed oriented simply-connected smooth $4$--manifold. We will describe the  branched covering in the statement of Theorem~\ref{thm1}, that is \linebreak $f\: T^2 \x \Sigma_g \to X$, as a composition of four simpler double branched coverings $f_1, f_2, f_3, f_4$.
While all the latter will be branched over embedded orientable surfaces, the branch locus of the composition will typically be singular.

For the clarity of the exposition, we will not explicitly keep track of how the topology is growing at each step, but instead, we will illustrate with some examples in Remark \ref{explicit} how one can deduce this information.

\medskip
\noindent \underline{\smash{\textsl{Step 1:}}} \, By Wall \cite{Wall}, the connected sum of $X$ with a certain number $m$ of copies of $S^2 \x S^2$ is diffeomorphic to a connected sum of copies of the standard $4$--manifolds $\CP$, $S^2 \x S^2$ and the $\K$ surface, taken with either orientations. Note that when $X$ is spin, the decomposition has only spin connected summands, and also that the resulting $4$--manifold does satisfy $11/8$ when $m$ is large enough. 

Moreover, since we have $\K \cs \Kb \cong \cs_{22} (S^2 \x S^2)$ and $\CP \cs (S^2 \x S^2) \cong \cs_2 \CP \cs \CPb$  \cite[Page 344]{GS}, the complete decomposition as above can be written as $\cs_a \K \cs_b (S^2 \x S^2)$ or $\cs_a \Kb \cs_b (S^2 \x S^2)$ when $X$ is spin (depending on the sign of the signature $\sigma(X)$), and as $\cs_a \CP \cs_b \CPb$ when $X$ is non-spin, for some non-negative integers $a$ and $b$ which are not both zero. In the spin case, we can guarantee that $b \geq 2a$ by taking sufficiently many stabilizations.

The conjugation map $(z_1, z_2) \mapsto (\bar z_1, \bar z_2)$, which is an anti-holomorphic involution on $\CPo \x \CPo \cong S^2 \x S^2$,  induces a double branched covering $S^2 \times S^2 \to S^4$, where the branch locus is the unknotted $T^2 \subset S^4$ (bounding a handlebody). Taking equivariant connected sum of $m$ copies of it, we get an involution on $\cs_m (S^2 \x S^2)$, which induces a double covering $\cs_m (S^2 \x S^2) \to S^4$  branched along an unknotted $\Sigma_m$, for every $m \geq 1$.

We can now take a double covering of $X$ branched along an unknotted $\Sigma_{m}$ in $X$  (viewing $X \cong X \cs S^4$, take  an unknotted $\Sigma_{m}$ in $S^4$), which we denote by $f_1\: X_1 \to X$, where clearly $X_1 \cong X \cs X\cs_m (S^2 \x S^2)$. We choose $m \geq 1$ such that $X\cs_m (S^2 \x S^2)$ completely decomposes, so does $X_1$ (as one gets at least $m$ copies of $S^2 \x S^2$ after decomposing  $X \cs_m (S^2 \x S^2)$). Then $X_1$ is diffeomorphic to one of the standard connected sums we listed above.

\smallskip
\noindent \underline{\smash{\textsl{Step 2:}}} \, We  would like to obtain a double branched  covering of $X_1$ by some $\cs_g (S^2 \x S^2)$. We will describe this covering in essentially two different ways, depending on whether $X$ (and thus $X_1)$ is spin or not.

The $\K$ surface can be obtained as a holomorphic double  covering of $S^2 \x S^2$ branched along a curve of bi-degree $(4,4)$ in $\CPo \x \CPo \cong S^2 \x S^2$ \cite[Page 262]{GS}. 
Reversing the orientations, we see that $\Kb$ is also a double branched covering of $S^2 \x S^2$ (recall that $S^2 \x S^2$ admits an orientation-reversing diffeomorphism). By taking equivariant connected sums, we can then express both $\cs_n \K$ and $\cs_n \Kb$ as  branched double coverings of $\cs_n (S^2 \x S^2)$. Taking \mbox{$n=2a$,} we then conclude that $\cs_a \K \cs_b (S^2 \x S^2)$ admits a double branched covering by  $\cs_{2a} \K \cs_{2a} \Kb \cs_{2(b-2a)} (S^2 \x S^2)$. Since $\K \cs \Kb \cong \cs_{22} (S^2 \x S^2)$, we have obtained the desired double branched cover $\cs_g (S^2 \x S^2)$, for $g=40a+2b$. Mirroring the same argument, we  see that $\cs_a \Kb \cs_b (S^2 \x S^2)$ is also double branched covered by some $\cs_g (S^2 \x S^2)$. This concludes the construction in the spin case.

The following variation can be run for both spin and non-spin manifolds.
Switching the two factors $(z_1, z_2) \mapsto (z_2, z_1)$, which is a holomorphic involution on $\CPo \x \CPo \cong S^2 \x S^2$,  induces a double branched covering $S^2 \times S^2 \to \CP$, where the branch locus is the quadric (this may be interpreted as the map taking a pair of numbers to the quadratic equation having those roots). Reversing the orientations, we obtain a double branched covering over $\CPb$. Taking equivariant connected sums once again, we then deduce that $\cs_a (S^2 \x S^2) \cs_b (S^2 \x S^2)$ is a double branched covering of $\cs_a \CP \cs_b \CPb$. So in the non-spin case, we arrive at the desired double covering $\cs_g (S^2 \x S^2)$ as well.

We let $f_2\: X_2 \to X_1$ denote the double branched covering we  described in either case.

\smallskip
\noindent \underline{\smash{\textsl{Step 3:}}} \, We next show that $S^1 \x \cs_g (S^1 \x S^2)$ is a double branched covering of $ \cs_g (S^2 \x S^2)$, which will prescribe our next covering $f_3 \: X_3 \to X_2$. A similar double branched covering over $\cs_g (\CP \cs \CPb)$ was described by Neofytidis in \cite[Theorem~1]{Neofytidis}. (Also see \cite{NeofytidisThesis} for similar constructions in other dimensions.)

The hyperelliptic involution on $T^2$  induces a double branched covering $p\: T^2 \to S^2$ with four simple branch points $x_1, x_2, x_3, x_4 \in S^2$. Taking its product with the identity map on $S^2$ yields  a double branched covering $p \x \textrm{id}_{S^2}  \: T^2 \x S^2 \to S^2 \x S^2$ with branch locus $\{x_1, x_2, x_3, x_4\} \x S^2$. Note that if $g=1$, we can stop here and skip Step 4.

Let $D \subset S^2$ be a $2$--disk containing exactly two branch points of $p$. So $A = p^{-1}(D) \subset T^2$ is an equivariant annulus that contains two fixed points of the hyperelliptic involution. Moreover, $D$ can be chosen such that $A$ is a union of fibers of the trivial $S^1$--bundle $T^2 = S^1 \times S^1 \to S^1$ given by the canonical projection onto the second factor.

Let $S^2 \cong S^2 \times \{y\} \subset S^2 \times S^2$ be a fiber sphere, for a certain $y \in S^2$. Let $D' \subset S^2$ be a disk centered at $y$. Then, $U = D \times D' \subset S^2 \times S^2$ is a fibered bidisk, whose preimage $V = (p\times \id_{S^2})^{-1} (U) \cong A \times D'$ is a fibered neighborhood of a fiber of the trivial $S^1$--bundle \[T^2 \x S^2 = S^1 \x (S^1 \x S^2) \to S^1 \x S^2.\]
By taking two copies of the branched covering $T^2 \times S^2 \to S^2 \times S^2$, and performing equivariant fiber sum upstairs along $V$ and connected sum downstairs along $U \cong D^4$, and repeating the construction for every $g \geq 2$, we finally get a branched double covering $S^1 \x \cs_g (S^1 \x S^2) \to \cs_g (S^2 \x S^2)$.

We can also describe this branched covering as follows: start with a double covering $q \: S^1 \times D^1 \to D^2$ branched over two points in $\Int D^2$ (this is the above branched covering $A \to D$), so the product $q \x \id_{D^1} \: S^1 \times D^1 \times D^1 \to D^2 \x D^1$ yields a double covering $q' \: S^1 \x D^2 \to  D^3$ branched over the union of two parallel proper trivial arcs in $D^3$ (this fills the above branched covering $p \: T^2 \to S^2$), up to the identifications $S^1 \times D^1 \times D^1 \cong S^1 \x D^2$ and $D^2 \x D^1 \cong D^3$. Then, we get a double branched covering \linebreak $q'' = q' \x \id_{S^2} \: S^1 \times D^2 \x S^2 \to D^3 \x S^2$.
Let $D \subset S^2$ be a $2$--disk. Up to the identification $D^2 \times D^1 \times S^2 \cong D^3 \times S^2$, we consider the bidisks $C^- = D^2 \x \{-1\} \x D$ and \mbox{$C^+ = D^2 \x \{1\} \x D \subset \Bd(D^3 \times S^2)$,} each of which intersects the branch locus of $q''$ along the union of two parallel proper trivial $2$--disks. Consider $g$ copies of $q''$, say $q''_i \: (S^1 \times D^2 \x S^2)_i \to (D^3 \x S^2)_i$, and let $C_i^-, C_i^+ \subset \Bd(D^3 \times S^2)_i$ be the corresponding bidisks. Thus, we obtain a double branched covering \[q''' = q''_1 \cup \dots\cup q''_g \: {\cup_i} (S^1 \times D^2 \x S^2)_i \to \cup_i (D^3 \x S^2)_i,\] where $(D^3 \x S^2)_i$ is attached to $(D^3 \x S^2)_{i+1}$ by identifying $C_i^+$ with $C_{i+1}^-$ and\linebreak $(S^1 \times D^2 \x S^2)_i$ is attached to $(S^1 \times D^2 \x S^2)_{i+1}$ by identifying $(q''_i)^{-1}(C_i^+)$ with\linebreak $(q''_{i+1})^{-1}(C_{i+1}^-)$ in the obvious way, for all $i = 1,\dots,g-1$. This in turn is a double branched covering \[q''' \: S^1 \times \csb_g(D^2 \x S^2) \to \csb_g (D^3 \x S^2),\] as it can be easily realized by looking at the attaching maps, where $\csb$ denotes the boundary connected sum. Finally, the desired branched covering $S^1 \x \cs_g (S^1 \x S^2) \to \cs_g (S^2 \x S^2)$ can be obtained by restricting $q'''$ to the boundary.

\smallskip
\noindent \underline{\smash{\textsl{Step 4:}}} \, Our final double branched covering $f_4\: T^2 \x \Sigma_g \to  S^1 \x  \cs_g (S^1 \x S^2)$ is a special case of Proposition~\ref{prop1} and can be  obtained by taking the product of the identity map on the $S^1$ factor with a double branched covering $S^1 \times \Sigma_g \to \cs_g (S^1 \x S^2)$. The latter can be  derived from the work of Sakuma \cite{Sakuma} we mentioned in the introduction, or from Montesinos' alternative construction \cite{Montesinos}, which is quicker to describe here: the involution $(z,t) \to (\bar{z}, -t)$ on the annulus $A=S^1\x [-1,1] \subset \C \x \R$ induces a double  covering $q\: A \to D^2$ branched at two points (this is same as the double branched cover described in Step 3), so we get a double branched covering $q \x \textrm{id}_{S^1}\: A \x S^1 \to D^2 \x S^1$. Then, for any open book decomposition of a closed connected oriented 3-manifold $Y$ with pages $\Sigma_k^m$ and monodromy $\phi$, we can get a double covering $h\: Y' \to Y$ branched over two parallel copies of the binding, where $Y'$ is now a surface bundle whose fiber and the monodromy are the doubles of $\Sigma_k^m$ and $\phi$. Indeed, by lifting the usual splitting $Y = (D^2 \x \Bd \Sigma_k^m) \cup_\Bd T(\phi)$ that gives the open book decomposition of $Y$, with the branch link contained in $D^2 \x \Bd \Sigma_k^m$, and where $T(\phi)$ denotes the mapping torus of $\phi$, one obtains a splitting $Y' = (A \x \Bd \Sigma_k^m) \cup_\Bd (T(\phi)_1\cup T(\phi)_2)$, with the annulus $A$ instead of $D^2$, where $T(\phi)_1$ and $T(\phi)_2$ are two disjoint copies of $T(\phi)$ (the branched covering $h \: Y' \to Y$ is trivial over $T(\phi)$). By looking at the attaching maps, it is immediate to get the bundle structure on $Y'$ as above.

In our case, since $\cs_g (S^1 \x S^2)$ admits a planar open book with pages $\Sigma_0^{g+1}$ and $\phi=\textrm{id}$, we obtain the desired covering. (The covering produced by the arguments of both Sakuma and Montesinos  in this  simple setting is equivalent to the one given in \cite[Proposition~4]{KotschickNeofytidis}.)

The composition $f=f_1 \circ f_2 \circ f_3 \circ f_4 : T^2 \x \Sigma_g \to X$ gives the desired covering.\\

\noindent \underline{\smash{\textsl{The spin case:}}} \, Let us conclude by observing that our construction is natural with respect to the spin structures, when $X$ is spin, and then briefly discuss the topology of the branch locus of $f$.

Recall that a spin structure on a $4$--manifold is the same as a trivialization of the tangent bundle over the $1$--skeleton that extends over the $2$--skeleton \cite{Milnor, Kirby}. We may use a handlebody decomposition in this definition. Given an unramified cover over a spin $4$--manifold the trivialization will lift to the tangent bundle of the cover restricted to the $1$--skeleton and any extension to the $2$--skeleton, so there is a natural lift of a spin structure to a covering space.

Now consider a $2$--fold branched covering with branch locus $B$. We may build a handle decomposition of the base in the following way. Start with a handle decomposition of $B$. This extends to a handle decomposition of a tubular neighborhood of $B$ with only zero, one and two handles. Now extend this to a handle decomposition of the rest of the base $X$. Finally turn the entire handle decomposition over. Notice that all of the $1$--handles of this new handle decomposition are in the exterior of $B$. Each of these handles lifts to the cover of the exterior and the restriction of the spin structure to the exterior lifts to the cover. We now complete the handle decomposition of the total space of the branched cover as follows. Use the identification of the inverse image $\widetilde B$ with $B$ to construct a decomposition of $\widetilde B$  which is then extended to a decomposition of the normal bundle of $\widetilde B$ . Turn this upside down and add it the to decomposition of the inverse image of the exterior. This only adds 2-, 3- and 4-handles to the decomposition. It is not necessarily true that the trivialization of the tangent bundle over the $1$--skeleton will extend over the $2$--skeleton. It will extend precisely when the mod two reduction of the integral homology class $[B]/2$ is zero in the second homology of the base with $\mathbb{Z}_2$ coefficients \cite{Brand, Nagami}. Note that the class of $[B]$ is necessarily divisible by 2 due to the existence of the double branched cover.  So, a spin structure does not have to lift to the total space of a $2$--fold branched covering, but if it does, there is a natural lift.

It is now straightforward to check that each double cover $f_i$ that we employed in our construction when $X$ is spin satisfies the above criterion, so for the initial spin structure $\mathfrak{s}$ on $X$, there is a spin structure $\mathfrak{s}'$ on $X' \cong T^2 \x \Sigma_g$ constructed this way. (Note that there are  $2^{2(g+1)}$ different spin structures on $X'$.) Thus the branched covering $X' \to X$  is  compatible with the spin structures $\mathfrak{s}$ on $X$ and $\mathfrak{s}'$ on $X'$.\\

\noindent \underline{\smash{\textsl{The branch locus:}}} \, The branch locus $B_f \subset X$ of $f$ is given by \[B_f = B_{f_1} \cup f_1\Big(B_{f_2} \cup f_2\big(B_{f_3} \cup f_3(B_{f_4})\big)\!\Big),\] where $B_{f_i} \subset X_{i-1}$ denotes the branch locus of $f_i$, for $i=1,2,3,4$, with $X_0 = X$. Each $B_{f_i}$ is a smooth embedded closed orientable surface in $X_{i-1}$. By taking into account that each covering $f_i$ is two-to-one and its tangent map has a 2-dimensional kernel along the branch set, an easy transversality argument based on perturbing the $f_i$'s up to isotopy, shows that the branch locus $B_f \subset X$ can be assumed to be a smooth orientable surface away from at most finitely many singular points, which are transversal or tangential double points (at the latter the local link has two trivial components with linking number $\pm 2$).
\qed

\bigskip
\section{Ancillary Remarks}

Let us list a few comments in relation to Theorem \ref{thm1}, its proof and related works.

\begin{rk}[Variations] \label{variations}
In Step 1 above we could have stabilized  by taking connected sums with copies of $\CP$ and  $\CPb$ so that we get a double covering \mbox{$g_1\: \#_a\CP \cs_b \CPb \to X$} branched over a genus $m$ \textit{non-orientable} surface, which is trivially embedded in $X$, for certain integers $a,b$ and $m$ (once again by Wall \cite{Wall}). The complex conjugation on $\CP$ induces this double covering $\CP \to S^4$ branched over the standard smooth $\RP \subset S^4$ \cite{Ma73, Ku74}.  
Now, we can invoke Theorem 1.2 in \cite{PiergalliniZuddas} to conclude that there exists a $4$--fold simple branched covering \,$g_2\:  \Sigma_h \times \Sigma_g \to \cs_a\CP \cs_b \CPb$ for \emph{every} given $a,b \geq 0$ and $h\geq 1$, and for some $g$ large enough. Thus, the composition $g_1 \circ g_2 \: \Sigma_g \times \Sigma_h \to X$ is a degree 8 branched covering.

Again by Theorem 1.2 in \cite{PiergalliniZuddas} (see also Remark 2 therein), there exist degree 4 branched coverings $T^4 = T^2 \times T^2 \to X$, with $X = \cs_m \CP \cs_n \CPb$ and $X = \cs_n (S^2 \x S^2)$, for every $m,n \leq 3$. Note that the case $X = S^2 \x S^2$ is straightforward by taking the product $p \x p \: T^2 \x T^2 \to S^2 \x S^2$, and the case $X = \cs_2 (S^2 \x S^2)$ was previously obtained by Rickman \cite{Rickman}. Branched coverings from the $n$--dimensional torus are relevant in connection with the theory of \textit{quasiregularly elliptic} manifolds, see Bonk and Heinonen \cite{Bonk}. In this direction, a result by Prywes \cite[Theorem 1.1]{Prywes} implies that if there is a branched covering $T^4 \to X$, then $b_1(X) \leq 4$ and $b_2(X) \leq 6$, so in Theorem \ref{thm1} we cannot take $g \leq 1$ if $b_2(X) \geq 7$.

However, unlike in our construction above, the results in \cite{PiergalliniZuddas} do not give explicit branched coverings, and there is not much control on the topology of the branch locus. 
\end{rk}

\smallskip
\begin{rk}[Branched cover geometries] \label{geometry}
Theorem~\ref{thm1} and our subsequent remark in the introduction imply that any $X$ with finite $\pi_1(X)$ is branched covered by $T^2 \x \Sigma_g$, where it is easy to see from our proof that we can always assume $g\geq 2$. In terms of \mbox{$4$--dimensional} geometries \cite{Hillman}, this shows that all such $X$ can be branched covered by a $4$--manifold with $\mathbb{E}^2 \x \mathbb{H}^2$ geometry. However, if we replace the double branched covering $h\: Y' \to \cs_g (S^1 \x S^2)$ we used in the construction of $f_4=\textrm{id}_{S^1} \x h$ with the one built by Brooks  in \cite{Brooks}, we can also get $Y'$ to be a $\Sigma_g$--bundle over $S^1$ with hyperbolic total space. Therefore, any $X$ with finite $\pi_1(X)$ can also be branched covered by a $4$--manifold with $\mathbb{E} \x \mathbb{H}^3$ geometry. Similarly, one can modify the construction in Proposition~\ref{prop1} to get a  double branched cover of any product $4$--manifold $S^1 \x Y$ by a $4$--manifold with $\mathbb{E} \x \mathbb{H}^3$ geometry.
\end{rk}

\smallskip
\begin{rk}[Topology of the branched coverings] \label{explicit}
Here we will try to demonstrate by way of example how one can control the topology of the branched coverings  in Theorem~\ref{thm1}. For some variety, we will run our construction for two infinite families of  irreducible $4$--manifolds which are not completely decomposable: Dolgachev surfaces, which are non-spin complex surfaces of general type, and knot surgered $\K$ surfaces of Fintushel--Stern, which include spin \mbox{$4$--manifolds} that do not admit any symplectic structures \cite{FintushelStern}. The members of either one of these two families of simply-connected $4$--manifolds completely decompose after a single stabilization by $S^2 \x S^2$; see e.g. \cite{Baykur}.  Now, if $X$ is a Dolgachev surface, we can  take the branch locus of $f_1$ as an unknotted $T^2$, and  get \mbox{$X_1= \cs_3 \CP \cs_{19} \CPb$.} The branched covering $f_2$ performed along the connected sum of $22$ quadrics (each coming from distinct copies of $\CP$ and $\CPb$) then gives $X_2=\cs_{22} (S^2 \x S^2)$, and the last two coverings yield $X'=T^2 \x \Sigma_{22}$. If we take $X$ to be a knot surgered $\K$ surface instead, thinking ahead of the second step, we take the branch locus of $f_1$ this time as $\Sigma_4$, so $X_1= \cs_2 \K \cs_4 (S^2 \x S^2)$. The next double cover $f_2$ is  taken along a connected sum of four bi-degree $(4,4)$ curves (each coming from distinct copies of $S^2 \x S^2$, with non-complex orientation), and we get $X_2= \cs_4 \K \cs_4 \Kb \cong \cs_{88} (S^2 \x S^2)$. The last two coverings this time yield $X'=T^2 \x \Sigma_{88}$.
\end{rk}

\smallskip
\begin{rk} [Symplectic domination]
\label{domination}
A recent article of Fine--Panov provides a symplectic domination result \cite[Theorem 1]{FinePanov} which is worth mentioning here. Their beautiful construction is very general: for any closed oriented even dimensional smooth manifold $M$, the authors build a closed symplectic manifold $S$ of the same dimension with a positive degree map $f\: S \to M$. In dimension $4$, where we can compare their result with ours in Theorem~\ref{thm1}, their symplectic manifold $S$ is constructed as a Donaldson hypersurface in the $6$--dimensional symplectic twistor space $Z$ of a negatively pinched manifold $N$, where the latter admits a degree one map \mbox{$g\: N\to M$.} The construction of $N$, with sectional curvature arbitrarily close to $-1$ is implicit, and relies on the recent works of Ontaneda involving rather intricate new techniques in Riemannian geometry. (The condition on the sectional curvature is to guarantee that the twistor space  $Z$ of $N$ is a symplectic $6$--manifold.) Secondly, the construction of a symplectic hypersurface $S$ in $N$, which is built through asymptotically holomorphic techniques  of \cite{Donaldson}, is also implicit and the smooth topology of $S$ is effectively impossible to control. Hence, one does not have any information on the smooth topology of the dominating symplectic $4$--manifold $S$, other than that it is of general type, i.e. of Kodaira dimension $2$ \cite{FinePanov}. Besides the very implicit nature of this construction, since the map $f\: S \to M$ factors through the degree one map $g$ above, the authors' domination is essentially never  a branched covering. Moreover, because the symplectic twistor space $Z$ is in fact known to be non-K\"ahler \cite{Reznikov}, the dominating symplectic $4$--manifold $S$ has a priori no reason to be a K\"ahler surface. On the other hand, the dominating symplectic $4$--manifold $X'=T^2 \x \Sigma_g$ of Theorem~\ref{thm1} is obviously a K\"{a}hler surface, and $X'$  in both Theorem~\ref{thm1} and Proposition~\ref{prop1} is of Kodaira dimension $-\infty$, $0$ or $1$, depending on whether this (possibly trivial) $\Sigma_g$--bundle over $T^2$, has fiber genus $g=0$, $1$ or $\geq 2$, respectively.

Domination is certainly distinct from branched covering as the following example shows. There is a degree one map from $\Sigma_4\times\Sigma_2$ to $\Sigma_3\times\Sigma_2$ given by the extension of the natural collapse of a copy of $\Sigma_1^1\times\Sigma_2$ to $\Sigma_0^1\times\Sigma_2$. However there can be no branched covering from $\Sigma_4\times\Sigma_2$ to $\Sigma_3\times\Sigma_2$ since the Gromov norm of the former is $24(4-1)(2-1)=72$, the Gromov norm of the latter is $24(3-1)(2-1)=48$ and the Gromov norm is super multiplicative with respect to degree \cite{BK}.
\end{rk}


\clearpage
\begin{thebibliography}{99999}

\bibitem{Baykur} R. I. Baykur, \textit{Dissolving knot surgered $4$-manifolds by classical cobordism arguments,} J. Knot Theory Ramifications, 27:5 (2018), 1871001.

\bibitem{BaykurSunukjian} R. I. Baykur and N. Sunukjian, \textit{Round handles, logarithmic transforms and smooth $4$-manifolds}, J. Topol. 6 (2013), no. 1, 49--63.

\bibitem{Bonk}
M. Bonk and J. Heinonen,
\textit{Quasiregular mappings and cohomology},
Acta Math. 186 (2001), no. 2, 219--238.

\bibitem{Brand} N. Brand, \textit{Necessary conditions for the existence of branched coverings,} Invent. Math. 54 (1979), 1--10.

\bibitem{Brooks} R. Brooks, \textit{On branched coverings of $3$-manifolds which fiber over the circle,} J. Reine Angew. Math. 362 (1985), 87--101.

\bibitem{BK} M. Bucher-Karlsson, \textit{The simplicial volume of closed manifolds covered by $\mathbb{H}^2\times\mathbb{H}^2$}, J. Topol. 1 (2008), no. 3, 584--602.

\bibitem{Church}
P.\, T. Church,
\textit{Differentiable open maps on manifolds}, 
Trans. Amer. Math. Soc. 109 (1963), 87--100. 

\bibitem{Donaldson} S. K. Donaldson, \textit{Symplectic submanifolds and almost-complex geometry}, J. Differential Geom. 44 (1996), no. 4, 666--705.

\bibitem{Eliashberg} Y. Eliashberg, \textit{Recent advances in symplectic flexibility}, Bull. Amer. Math. Soc. (N.S.) 52 (2015), no. 1, 1--26.

\bibitem{FinePanov} J. Fine and D. Panov, \textit{Symplectic domination}, 
Bull. London. Math. Soc. doi:10.1112/blms.12402, 2020.

\bibitem{FintushelStern} R. Fintushel and R. Stern, \textit{Knots, links, and $4$-manifolds}, Invent. Math. 134 (1998), no. 2, 363--400. 

\bibitem{GS} R. Gompf and A. Stipsicz, \textit{$4$--manifolds and Kirby calculus}, Graduate Studies in Mathematics, vol. 20, American Mathematical Society, Providence, RI. 1999.

\bibitem{Hillman} J. Hillman, \textit{Four-manifolds, geometries and knots}, Geom. Topol. Monographs, 5 (2002), 379 pp. 

\bibitem{Kirby} R. Kirby, \textit{The topology of 4-manifolds}, Lecture Notes in Mathematics, 1374. (1989), 108 pp.

\bibitem{KirbyList} R. Kirby,  \textit{Problems in Low-Dimensional Topology}  (1995),  \url{https://www.math.berkeley.edu/~kirby/}. 

\bibitem{KotschickLoh} D. Kotschick and C. L\"{o}h, \textit{Fundamental classes not representable by products,} J. London Math. Soc. (2) 79 (2009), 545--561.

\bibitem{KotschickNeofytidis} D. Kotschick and C. Neofytidis, \textit{On three-manifolds dominated by circle bundles,} Math. Z. 274 (2013) 21--32.

\bibitem{Ku74}
N.\,H. Kuiper,
\textit{The quotient space of $\CP$ by complex conjugation is the $4$-sphere},
Math. Ann. 208 (1974), 175--177.

\bibitem{Ma73}
W.\, S. Massey,
\textit{The quotient space of the complex projective plane under conjugation is a 4-sphere}, 
Geometriae Dedicata 2 (1973), 371--374.

\bibitem{Milnor}
J. Milnor, \textit{Spin structures on manifolds}, Enseign. Math. 9 (1963), 198--203. 

\bibitem{Montesinos} J. M. Montesinos, \textit{On $3$-manifolds having surface bundles as branched coverings}, Proc. Amer. Math. Soc. 101 (1987), no. 3, 555--558.

\bibitem{Nagami} S. Nagami, \textit{On spin structures of double branched covering spaces}, JP J. Geom. Topol. 14 (2013), no. 2, 119--147.

\bibitem{Neofytidis} C. Neofytidis, \textit{Branched coverings of simply-connected manifolds}, Topology Appl. 178 (2014), 360--371.

\bibitem{NeofytidisThesis} C. Neofytidis, \textit{Non-zero degree maps between manifolds and groups presentable by products}, Dissertation (2014), LMU M\"{u}nchen. 

\bibitem{PiergalliniZuddas} R. Piergallini and D. Zuddas, \textit{Branched coverings of $\CP$ and other basic $4$-manifolds}, to appear in Bull. London Math. Soc., \url{https://doi.org/10.1112/blms.12463}.

\bibitem{Prywes}
E. Prywes,
\textit{A bound on the cohomology of quasiregularly elliptic manifolds},
Ann. of Math. 189 (2019), 863--883.

\bibitem{Reznikov} A. G. Reznikov, \textit{Symplectic twistor spaces}, Ann. Global Anal. Geom. 11 (1993), no. 2, 109--118. 

\bibitem{Rickman}
S. Rickman,
\textit{Simply connected quasiregularly elliptic 4-manifolds},
Ann. Acad. Sci. Fenn. Math. 31 (2006), no. 1, 97--110. 

\bibitem{Sakuma} M. Sakuma, \textit{Surface bundles over $S^1$ which are $2$-fold branched cyclic coverings of $S^3$}, Math. Sem. Notes 9 (1981), 159--180.

\bibitem{Thurston} W. P. Thurston, \textit{Some simple examples of symplectic manifolds},
Proc. Amer. Math. Soc. 55 (1976), no.2, 467--468.

\bibitem{Wall} C. T. C. Wall, \textit{On simply-connected $4$-manifolds,} J. London Math. Soc. 39 (1964) 141--149.

\end{thebibliography}
\end{document}